\documentclass[11pt, leqno]{article}
\usepackage{amsfonts}
\usepackage{latexsym}
\usepackage{amsmath}
\usepackage{amssymb}
\usepackage{amsthm}
\newtheorem{theorem}{Theorem}[section]
\newtheorem{lemma}[theorem]{Lemma}

\newtheorem{proposition}[theorem]{Proposition}

\newtheorem{remark}[theorem]{Remark}

\theoremstyle{definition}
\newtheorem{definition}[theorem]{Definition}

\theoremstyle{remark}

\newtheorem*{note*}{Note}

\numberwithin{equation}{section}

\makeatletter
\newcommand{\ls}{\leqslant}
\newcommand{\gr}{\geqslant}
\makeatother

\begin{document}
\small

\title{\bf A remark on the slicing problem}

\author{A. Giannopoulos, G. Paouris and B-H. Vritsiou}

\date{}

\maketitle

\begin{abstract}
\footnotesize The purpose of this article is to describe a reduction
of the slicing problem to the study of the parameter
$I_1(K,Z_q^{\circ }(K))=\int_K\|\langle \cdot ,x\rangle
\|_{L_q(K)}dx$. We show that an upper bound of the form
$I_1(K,Z_q^{\circ }(K))\ls C_1q^s\sqrt{n}L_K^2$, with $1/2\ls s\ls
1$, leads to the estimate $$L_n\ls \frac{C_2\sqrt[4]{n}\log\!
n}{q^{\frac{1-s}{2}}},$$ where $L_n:=\max\{ L_K:K\;\hbox{is an
isotropic convex body in}\,{\mathbb R}^n\}$.
\end{abstract}

\section{Introduction}

A convex body $K$ in ${\mathbb R}^n$ is called isotropic if it has
volume $|K|=1$, it is {\it centered}, i.e. its center of mass is
at the origin, and if its inertia matrix is a multiple of the
identity. The last property is equivalent to the existence of a
constant $L_K >0$ such that
\begin{equation}\int_K\langle x,\theta\rangle^2dx =L_K^2\end{equation}
for every $\theta $ in the Euclidean unit sphere $S^{n-1}$. It is
not hard to see that for every convex body $K$ in ${\mathbb R}^n$
there exists an affine transformation $T$ of ${\mathbb R}^n$ such
that $T(K)$ is isotropic. Moreover, this isotropic image is unique
up to orthogonal transformations; consequently, one may define the
isotropic constant $L_K$ as an invariant of the affine class of $K$.

The isotropic constant is closely related to the hyperplane
conjecture (also known as the slicing problem) which asks if there
exists an absolute constant $c>0$ such that
$\max\limits_{\theta\in S^{n-1}}|K\cap\theta^{\perp }|\gr c$ for every convex body $K$ of
volume $1$ in ${\mathbb R}^n$ with center of mass at the origin.
This is because, by Brunn's principle, for any convex body $K$ in
${\mathbb R}^n$ and any $\theta\in S^{n-1}$, the function $t\mapsto
|K\cap (\theta^{\perp }+t\theta )|^{\frac{1}{n-1}}$ is concave on
its support, and this is enough to show that
\begin{equation}\int_K\langle x,\theta\rangle^2dx\simeq |K\cap \theta^{\perp }|^{-2}.\end{equation}
Using this relation one can check that an affirmative answer to the slicing
problem is equivalent to the following statement: ``There exists an
absolute constant $C>0$ such that $L_K\ls C$ for every convex body
$K$''. We refer to the article \cite{MP} of Milman and Pajor for
background information about isotropic convex bodies.

It is known that $L_K\gr L_{B_2^n}\gr c>0$ for every convex body $K$
in ${\mathbb R}^n$ (we use the letters $c,c_1,C$ etc. to denote
absolute constants). In the opposite direction, let us write $L_n$ for the maximum of all isotropic
constants of convex bodies in ${\mathbb R}^n$,
\begin{equation}L_n:=\max\{ L_K:K\ \hbox{is isotropic in}\ {\mathbb R}^n\}.\end{equation}
Bourgain first proved in \cite{Bou} that
$L_n\ls c\sqrt[4]{n}\log\! n$ and, a few years ago, Klartag
\cite{Kl} obtained the estimate $L_n\ls c\sqrt[4]{n}$
(see also \cite{KlEM} for a second proof of this bound).

The purpose of this article is to describe a reduction of the
slicing problem (or, equivalently, the question whether $L_n$ can be bounded
by a quantity independent of the dimension $n$), to the study of the parameter
\begin{equation}I_1(K,Z_q^{\circ}(K))=\int_K\|\langle \cdot ,x\rangle \|_{L_q(K)}dx\end{equation}
for isotropic convex bodies $K$. Generally, if $K$ is a centered convex
body of volume $1$ in ${\mathbb R}^n$, then for every symmetric
convex body $C$ in ${\mathbb R}^n$ and for every $q\in (-n,\infty
)$, $q\neq 0$, we define
\begin{equation}I_q(K,C):=\left( \int_K\|x\|_C^qdx\right)^{1/q}.\end{equation}
The notation $I_1(K,Z_q^{\circ}(K))$ is then justified by the fact
that $\|\langle \cdot ,x\rangle\|_{L_q(K)}$ is the norm induced on
${\mathbb R}^n$ by the polar body $Z_q^{\circ }(K)$ of the
$L_q$-centroid body of $K$ (see the next section for background
information on $L_q$-centroid bodies).

Our reduction can be viewed as a continuation of Bourgain's approach
to the slicing problem in \cite{Bou}: the bound $O(\sqrt[4]{n}\:\!\log\! n)$
followed from the inequality
\begin{equation}nL_K^2\ls I_1\bigl(K,(T(K))^{\circ }\bigr),\end{equation}
after obtaining an upper bound for the quantity $I_1\bigl(K,(T(K))^{\circ}\bigr)$,
where $T\in SL(n)$ is a symmetric, positive definite matrix
such that the mean width of $T(K)$ satisfies the estimate
$w(T(K))\ls c\sqrt{n}\log\! n$ (the existence of such a
position for $K$ is guaranteed by Pisier's estimate on the norm of
the Rademacher projection; see \cite{Pi}). In Section 4 we prove the
following statement:

\begin{theorem}There exists an absolute constant $\rho\in (0,1)$
with the following property: given $\kappa ,\tau\gr 1$, for every
$n\gr n_0(\tau )$ and every isotropic convex body $K$ in ${\mathbb R}^n$
which satisfies the following entropy estimate:
\begin{equation}\label{eq:Th1.1_1}
\log N(K,tB_2^n)\ls \frac{\kappa n^2\log^2\! n}{t^2}\ \, \hbox{for all}\  t\gr \tau\sqrt{n\log\! n},
\end{equation}
we have that, if $q\gr 2$ satisfies
\begin{equation}\label{eq:Th1.1_2} 2\ls q\ls \rho^2n\  \hbox{and}\ I_1(K,Z_q^{\circ }(K))\ls \rho nL_K^2,\end{equation}
then
\begin{equation}
L_K^2\ls C\kappa\sqrt{\frac{n}{q}}\log^2\! n\,\max\left\{1,\frac{I_1(K,Z_q^{\circ}(K))}{\sqrt{qn}L_K^2}\right\}.
\end{equation}
\end{theorem}

Theorem 1.1 can lead to an upper bound for $L_n$, provided that
there exist {\it $(\kappa ,\tau )$-regular} isotropic convex bodies
in ${\mathbb R}^n$, i.e. bodies which satisfy the entropy estimate
$(\ref{eq:Th1.1_1})$ for a pair of constants $\kappa, \tau$, and at the same time
have maximal isotropic constant, i.e. $L_K\simeq L_n$. The
existence of such bodies is essentially established by \cite[Theorem 5.7]{DP}.
In Section 5 we give a self-contained proof of this fact; see Theorem 5.1.

Observe that, for every isotropic convex body $K$ in ${\mathbb R}^n$,
we have that both conditions in $(\ref{eq:Th1.1_2})$ are satisfied with $q=2$,
since $I_1(K,Z_2^{\circ}(K))\ls\sqrt{n}L_K^2$. Therefore,
Theorem 1.1 will give us that
\begin{equation}\label{eq:bound with q=2}
L_K^2\ls C_1\sqrt{n}\log^2n
\end{equation}
for any such body which is regular. Theorem 5.1 then guarantees that, for
some absolute constants $\kappa ,\tau $ and $\delta >0$, there
exists a $(\kappa ,\tau)$-regular isotropic convex body $K$ in
${\mathbb R}^n$ with $L_K\gr \delta L_n$, and hence $(\ref{eq:bound with q=2})$
leads us to Bourgain's bound again: $L_n\ls C_2\sqrt[4]{n}\log\! n$.

\smallskip

However, the behaviour of $I_1(K,Z_q^{\circ }(K))$ may allow us to
use much larger values of $q$. In Section 3 we discuss upper and
lower bounds for this quantity. For every isotropic convex body $K$
in ${\mathbb R}^n$ we have some simple general estimates:
\begin{enumerate}
\item[(i)]
For every $2\ls q\ls n$,
$$c_1\max\left\{ \sqrt{n}L_K^2, \sqrt{qn}, R(Z_q(K))L_K\right\} \ls I_1(K,Z_q^{\circ }(K))\ls c_2q\sqrt{n}L_K^2.$$
\item[(ii)]
If $2\ls q\ls \sqrt{n}$, then
$$c_1\max\left\{ \sqrt{n}L_K^2, \sqrt{qn}L_K\right\} \ls I_1(K,Z_q^{\circ }(K))\ls c_2q\sqrt{n}L_K^2.$$
\end{enumerate}
Any improvement of the exponent of $q$ in the upper bound
$I_1(K,Z_q^{\circ }(K))\ls cq\sqrt{n}L_K^2$ would lead to an
estimate $L_n\ls Cn^{\alpha}$ with $\alpha <\frac{1}{4}$. It seems
plausible that one could even have $I_1(K,Z_q^{\circ }(K))\ls c\sqrt{qn}L_K^2$,
at least when $q$ is small, say $2\ls q\ll \sqrt{n}$. Some evidence
is given by the following facts:
\begin{enumerate}
\item[(iii)]
If $K$ is an unconditional isotropic convex body in ${\mathbb R}^n$, then
$$c_1\sqrt{qn}\ls I_1(K,Z_q^{\circ }(K))\ls c_2\sqrt{qn}\log\! n$$
for all $2\ls q\ls n$.
\item[(iv)]
If $K$ is an isotropic convex body in ${\mathbb R}^n$ then, for every $2\ls q\ls\sqrt{n}$,
there exists a set $A_q\subseteq O(n)$ with $\nu (A_q)\gr 1-e^{-q}$ such that
$I_1(K,Z_q^{\circ}(U(K)))\ls c_3\sqrt{qn}\,L_K^2$ for all $U\in A_q$.
\end{enumerate}
\noindent The proofs of (i)-(iv) are given in Section 3.

\smallskip

We can make a final observation about the reduction of Theorem 1.1 on the basis that
there exist $(\kappa ,\tau)$-regular isotropic convex bodies $K$
in ${\mathbb R}^n$ with $L_K\gr \delta L_n$ (where
$\kappa,\tau,\delta >0$ are absolute constants) which, at the same
time, have ``small diameter": they satisfy $K\subseteq
\gamma\sqrt{n}L_K\,B_2^n$, where $\gamma >0$ is an absolute
constant (see Theorem 5.9). In Section 6, we show that then
it is enough to study the parameter $I_1(K,Z_q^{\circ}(K))$ within the class ${\cal
IK}_{{\rm sd}}$ of isotropic convex bodies which are
$O(\gamma)$-close to the Euclidean ball $D_n$ of volume $1$ and
have uniformly bounded isotropic constant. The precise statement which we prove
is the following: if we have an isotropic symmetric convex body $K$ in
${\mathbb R}^n$ satisfying $K\subseteq \gamma\sqrt{n}L_K\,B_2^n$, then
we can find an isotropic symmetric convex body $C$ such that
$L_C\ls c_1$, $c_2D_n \subseteq C \subseteq c_3\gamma D_n$, and
\begin{equation}\frac{I_1(K,Z_q^{\circ}(K))}{\sqrt{qn}L_K^2} \ls c_4\frac{I_1(C,Z_q^{\circ}(C))}{\sqrt{qn}}\end{equation}
for all $1\ls q\ls n$, where $c_1,c_2,c_3, c_4>0$ are absolute constants.

\section{Notation and preliminaries}

We work in ${\mathbb R}^n$, which is equipped with a Euclidean
structure $\langle\cdot ,\cdot\rangle $. We denote the corresponding
Euclidean norm by $\|\cdot \|_2$, and write $B_2^n$ for the Euclidean unit ball, and
$S^{n-1}$ for the unit sphere. Volume is denoted by $|\cdot |$. We
write $\omega_n$ for the volume of $B_2^n$ and $\sigma $ for the
rotationally invariant probability measure on $S^{n-1}$. We also
denote the Haar measure on $O(n)$ by $\nu $. The Grassmann manifold
$G_{n,k}$ of $k$-dimensional subspaces of ${\mathbb R}^n$ is
equipped with the Haar probability measure $\mu_{n,k}$. Let $k\ls n$
and $F\in G_{n,k}$. We will denote the orthogonal
projection from $\mathbb R^{n}$ onto $F$ by $P_F$. We also define
$B_F:=B_2^n\cap F$ and $S_F:=S^{n-1}\cap F$.

The letters $c,c^{\prime }, c_1, c_2$ etc. denote absolute positive
constants whose value may change from line to line. Whenever we write
$a\simeq b$, we mean that there exist absolute constants $c_1,c_2>0$
such that $c_1a\ls b\ls c_2a$.  Also if $K,L\subseteq \mathbb R^n$
we will write $K\simeq L$ if there exist absolute constants $c_1,
c_2>0$ such that $c_{1}K\subseteq L \subseteq c_{2}K$.

\medskip

A convex body in ${\mathbb R}^n$ is a compact convex subset $C$ of
${\mathbb R}^n$ with nonempty interior. We say that $C$ is
symmetric if $x\in C$ implies that $-x\in C$. We say that $C$ is
centered if it has center of mass at the origin, i.e. $\int_C\langle
x,\theta\rangle \,d x=0$ for every $\theta\in S^{n-1}$. The support
function of a convex body $C$ is defined by
\begin{equation}h_C(y):=\max \{\langle x,y\rangle :x\in C\},\end{equation}
and the mean width of $C$ is
\begin{equation}w(C):=\int_{S^{n-1}}h_C(\theta )\sigma (d\theta ).\end{equation}
For each $-\infty < q<\infty $, $q\neq 0$, we define the $q$-mean
width of $C$ by
\begin{equation}w_q(C):=\left(\int_{S^{n-1}}h_{C}^q(\theta)\sigma (d\theta )\right)^{1/q}.\end{equation}
The radius of $C$ is the quantity $R(C):=\max\{ \| x\|_2:x\in C\}$. Also, if the origin is
an interior point of $C$, the polar body $C^{\circ }$ of $C$ is defined as follows:
\begin{equation}C^{\circ}:=\{y\in {\mathbb R}^n: \langle x,y\rangle \ls 1\;\hbox{for all}\; x\in C\}.\end{equation}
Finally, we write $\overline{C}$ for the homothetic image of volume $1$
of a convex body $C\subseteq \mathbb R^n$, i.e. $\overline{C}:=\frac{C}{|C|^{1/n}}$.

\smallskip

Recall that if $A$ and $B$ are nonempty sets in $\mathbb R^n$, then
the covering number $N(A,B)$ of $A$ by $B$ is defined to be the
smallest number of translates of $B$ whose union covers $A$. In
this paper, $B$ will usually be a multiple of the Euclidean ball:
in those cases we also require that the centres of the translates
of $B$ are taken from the set $A$; one can easily check that this
additional requirement does not crucially affect our entropy estimates.

\medskip

\subsection{$L_q$-centroid bodies}

Let $K$ be a convex body of volume $1$ in ${\mathbb R}^n$. For every
$q\gr 1$ and every $y \in {\mathbb R}^n$ we set
\begin{equation}h_{Z_q(K)}(y):= \left(\int_K |\langle x,y\rangle|^{q}dx \right)^{1/q}.\end{equation}
The $L_q$-centroid body $Z_q(K)$ of $K$ is the
centrally symmetric convex body with support function
$h_{Z_{q}(K)}$. Note that $K$ is isotropic if and only if it is centered and $Z_{2}(K)=
L_{K}B_2^n$. Also, if $T\in GL(n)$ with ${\rm det}\,T = \pm 1$,
then $Z_{p}(T(K))= T(Z_{p}(K))$. From H\"{o}lder's inequality it follows that
$Z_1(K)\subseteq Z_p(K)\subseteq Z_q(K)\subseteq Z_{\infty }(K)$ for
all $1\ls p\ls q\ls \infty $, where $Z_{\infty }(K)={\rm conv}\{K,-K\}$.
Using Borell's lemma (see \cite[Appendix III]{MS}),
one can check that inverse inclusions also hold:
\begin{equation}\label{eq:Zq_1} Z_q(K)\subseteq \beta_1 qZ_1(K),\end{equation}
and more generally,
\begin{equation}\label{eq:Zq_2} Z_q(K)\subseteq \beta_2\frac{q}{p}Z_p(K)\end{equation}
for all $1\ls p<q$. In particular, if $K$ is isotropic, then
$R(Z_q(K))\ls \beta_1qL_K$. One can also check that if $K$ is
centered, then $Z_q(K)\supseteq \beta_3\,K$ for all $q\gr n$ (see
\cite{Pa1}). All the constants $\beta_i,\overline{\beta }_j$ that
appear in this section are absolute positive constants which may
be used again in the arguments of the next sections.

\medskip

Let $C$ be a symmetric convex body in ${\mathbb R}^n$ and let
$\|\cdot\|_C$ denote the norm induced on ${\mathbb R}^n$ by $C$. The
parameter $k_{\ast }(C)$ is the largest positive integer $k\ls n$
with the property that the measure of $F\in G_{n,k}$ for which we have $\frac{1}{2}
w(C)B_F \subseteq P_F(C) \subseteq 2 w(C)B_F$ is greater than
$\frac{n}{n+k}$. It is known that
\begin{equation}\beta_4n\frac{w(C)^2}{R(C)^2} \ls  k_{\ast }(C)\ls  \beta_5n \frac{w(C)^2}{R(C)^2}.\end{equation}
The $q$-mean width $w_q(C)$ is equivalent to $w(C)$ as long as $q\ls k_{\ast}(C)$.
Litvak, Milman and Schechtman proved in \cite{LMS} that, for
every symmetric convex body $C$ in $\mathbb R^n$,
\begin{enumerate}
\item[(i)]
If $1\ls  q\ls  k_{\ast }(C)$ then $w(C) \ls w_q(C) \ls \beta_6w(C)$.
\item[(ii)]
If $k_{\ast }(C)\ls  q\ls n$ then $\beta_7\sqrt{q/n}\, R(C) \ls  w_q(C) \ls  \beta_8 \sqrt{q/n}\, R(C)$.
\end{enumerate}

\smallskip

Let $K$ be a centered convex body of volume $1$ in ${\mathbb
R}^n$. Recall that, for every symmetric convex body $C$ in
${\mathbb R}^n$ and for every $q\in (-n,\infty )$, $q\neq 0$, we
define
\begin{equation}I_q(K,C):=\left( \int_K\|x\|_C^qdx\right)^{1/q}.\end{equation}
When $C=B_2^n$, we write $I_q(K):=I_q(K,B_2^n)$ for simplicity. In
\cite{Pa3} and \cite{Pa4} it is proved that for every $1\ls q\ls
n/2$,
\begin{equation}\label{eq:Zq_3}
I_q(K) \simeq \sqrt{n/q}\,w_q(Z_q(K))\ \hbox{and}\  I_{-q}(K) \simeq \sqrt{n/q}\, w_{-q}(Z_q(K)).
\end{equation}
The parameter $q_{\ast}(K)$ is also defined by
\begin{equation}q_{\ast}(K):= \max\{ q\ls n : k_{\ast}(Z_q(K)) \gr q\}.\end{equation}
Then, the main result of \cite{Pa4} states that,
for every centered convex body $K$ of volume $1$ in ${\mathbb R}^n$,
one has $I_{-q}(K) \simeq I_q(K)$ for every $1\ls q \ls q_{\ast}(K)$.
In particular, for all $q\ls q_{\ast}(K)$ one has
$I_q(K) \ls \beta_9I_2(K)$. If $K$ is isotropic, one can check that
$q_{\ast }(K)\gr c\sqrt{n}$, where $c>0$ is an absolute constant
(for a proof, see \cite{Pa3}). Therefore,
\begin{equation}\label{eq:Zq_4}
I_q(K)\ls \beta_{10}\sqrt{n}L_K\;\;\hbox{for every}\;q\ls\sqrt{n}.
\end{equation}
In particular, from $(\ref{eq:Zq_3})$ and $(\ref{eq:Zq_4})$ we see that $w(Z_q(K))\simeq
w_q(Z_q(K))\simeq\sqrt{q}L_K$ for all $q\ls\sqrt{n}$.

\medskip

\subsection{The bodies $B_q(K,F)$}

Another family of convex bodies associated with a centered
convex body $K \subset {\mathbb R}^n$ was introduced by Ball in
\cite{Ball} (see also \cite{MP}): to define them, let us consider a
$k$-dimensional subspace $F$ of ${\mathbb R}^n$ and its orthogonal
subspace $E$. For every $\phi\in F \setminus \{0\}$ we denote by $E^+(\phi)$
the halfspace $\{x\in {\rm span}\{E, \phi \}:\langle x,\phi\rangle\gr 0\}$.
Ball proved that, for every $q\gr 0$, the function
\begin{equation}
\phi\mapsto \| \phi\|_2^{1+\frac{q}{q+1}}\left( \int_{K \cap E^+(\phi)} \langle x,\phi \rangle^q dx\right)^{-\frac{1}{q+1}}
\end{equation}
is the gauge function of a convex body $B_q(K,F)$ on $F$. Several
properties of these bodies can be found in \cite{Ball}, \cite{MP}
and also in \cite{Pa3}, \cite{Pa4}. In Section 5, we will make use
of only two of those:
\begin{enumerate}
\item[(i)]
Let $K \subset {\mathbb R}^n$ be isotropic, let $1\ls k <n$ and let $F\in G_{n,k}$.
Then the body $\overline{B}_{k+1}(K,F)$ is {\it almost isotropic}, namely it has
(by definition) volume 1, and we can write
$\overline{B}_{k+1}(K,F)\simeq T(\overline{B}_{k+1}(K,F))$
where $T(\overline{B}_{k+1}(K,F))$ is an isotropic (in the regular sense)
linear image of $\overline{B}_{k+1}(K,F)$. In addition,
\begin{equation}\label{eq:KB1}
|K\cap F^{\perp}|^{1/k} \simeq \frac{L_{B_{k+1}(K,F)}}{L_K}.
\end{equation}
\item[(ii)]
Let $K, F$ and $k<n$ be as above and consider any $p\in [1,k]$. Then
\begin{equation}\label{eq:KB2}
Z_p \bigl(\overline{B}_{k+1}(K,F)\bigr)\simeq |K\cap F^{\perp}|^{1/k} P_F(Z_p(K)).
\end{equation}
\end{enumerate}

\medskip

\subsection{Two related lemmas}

We close this section with two lemmas that will be used in the
sequel; they reveal some properties of the support function of the
$L_q$-centroid bodies of a convex body with respect to subsets
or certain integrals of maxima.

\medskip

\begin{lemma}Let $K$ be a convex body of volume $1$ in ${\mathbb R}^n$,
and consider any points $z_1,z_2,\ldots ,z_N\in \mathbb R^n$. If $q\gr 1$
and $p\gr \max \{\log\! N, q\}$, then
\begin{equation}
\left(\int_K\max_{1\ls i\ls N} |\langle x,z_i\rangle |^qdx\right)^{1/q} \ls \overline{\beta }_1\max_{1\ls i\ls N}h_{Z_p(K)}(z_i),
\end{equation}
where $\overline{\beta }_1>0$ is an absolute constant.
\end{lemma}

\noindent {\it Proof.} Let $p\gr \max \{\log\! N, q\}$ and $\theta\in S^{n-1}$.
Markov's inequality shows that
\begin{equation}
|\{ x\in K: |\langle x,\theta\rangle |\gr e^3h_{Z_p(K)}(\theta )\}|\ls e^{-3p}.
\end{equation}
Since $x\mapsto |\langle x,\theta\rangle |$ is a seminorm,
from Borell's lemma (see \cite[Appendix III]{MS}) we get that
\begin{equation}
|\{x\in K: |\langle x,\theta\rangle |\gr e^3th_{Z_p(K)}(\theta )\}|\ls (1-e^{-3p})\left(\frac{e^{-3p}}{1-e^{-3p}}\right)^{\frac{t+1}{2}}\ls e^{-pt}
\end{equation}
for every $t\gr 1$. We set $S:=e^3\max\limits_{1\ls i\ls
N}h_{Z_p(K)}(z_i)$. Then, for every $t\gr 1$ we have that
\begin{align}
|\{x\in K:\max_{1\ls i\ls N}|\langle x,z_i\rangle |\gr St\}| 
&\ls \sum_{i=1}^N|\{x\in K: |\langle x,z_i\rangle |\gr e^3th_{Z_p(K)}(z_i)\}|
\\ \nonumber
&\ls Ne^{-pt}.
\end{align}
\noindent
It follows that
\begin{align}
\int_K\max_{1\ls i\ls N}|\langle x,z_i\rangle |^qdx 
&= q\int_{0}^{\infty}s^{q-1}|\{x\in K:\max_{1\ls i\ls N} |\langle x, z_i\rangle |\gr s\}| \,ds
\\ \nonumber
&\ls S^q + q\int_S^{\infty} s^{q-1}|\{x\in K:\max_{1\ls i\ls N} |\langle x,z_i\rangle |\gr s\}| \, ds
\\ \nonumber
&=   S^{q}\left( 1+ q\int_{1}^{\infty}t^{q-1}|\{x\in K:\max_{1\ls
i\ls N}|\langle x, z_{i}\rangle |\gr St\}|\,dt\right)
\\ \nonumber
&\ls S^q\left(1+ qN \int_{1}^{\infty} t^{q-1}e^{-pt} dt\right)
\\ \nonumber
&=   S^q\left(1+ \frac{qN}{p^q} \int_{p}^{\infty} t^{q-1}e^{-t} dt\right)
\\ \nonumber
&\ls S^q\left( 1+\frac{qN}{p^q} e^{-p} p^q \right)
\\ \nonumber
&\ls (3S)^{q},
\end{align}
\noindent
where we have also used the fact that, for every $p\gr q\gr 1$,
\begin{equation}\label{eqp:GammaGrowth}\int_{p}^{\infty} t^{q-1}e^{-t}dt \ls e^{-p} p^{q}.\end{equation}
This finishes the proof (with $\overline{\beta }_1=3e^3$).
$\hfill\Box $

\smallskip

\begin{remark}
{\rm It is a well-known fact (see e.g. \cite[Proposition 2.5.1]{Gian}) that
\begin{equation}
\int_K \max_{1\ls i\ls N}|\langle x,z_i \rangle| dx \ls C_1\log\! N \max_{1\ls i\ls N}h_{Z_1(K)}(z_i).
\end{equation}
Through a variant of the argument in \cite{Gian}, and using $(\ref{eqp:GammaGrowth})$ as well, one can also show that for $q\ls \log\! N$,
\begin{equation}
\left(\int_K \max_{1\ls i\ls N}|\langle x,z_i \rangle|^q dx\right)^{1/q} \ls C_2\log\! N \max_{1\ls i\ls N}h_{Z_1(K)}(z_i).
\end{equation}
Now, both inequalities can be directly deduced from Lemma 2.1 combined with $(\ref{eq:Zq_1})$, however the lemma provides additional information on how well the quantities $\left(\int_K \max_i |\langle x,z_i \rangle|^q dx\right)^{1/q}$
and $\max_i\left(\int_K |\langle x,z_i \rangle|^q dx\right)^{1/q} \equiv \max_i h_{Z_q(K)}(z_i)$ can be compared: for $q\ls \log\! N$, using also $(\ref{eq:Zq_2})$, we have that
\begin{equation}
\left(\int_K \max_{1\ls i\ls N}|\langle x,z_i \rangle|^q dx\right)^{1/q} \ls \overline{\beta}_1\max_{1\ls i\ls N}h_{Z_{\log\! N}(K)}(z_i)\ls C\frac{\log\! N}{q} \max_{1\ls i\ls N}h_{Z_q(K)}(z_i),
\end{equation}
whereas for $q > \log\! N$,
\begin{equation}
\left(\int_K\max_{1\ls i\ls N} |\langle x,z_i\rangle |^qdx\right)^{1/q}\simeq \max_{1\ls i\ls N}\left(\int_K |\langle x,z_i\rangle |^qdx\right)^{1/q}.
\end{equation}}
\end{remark}

\medskip

We now turn our attention to the $L_q$-centroid bodies of subsets of $K$.

\begin{lemma}
Let $K$ be a convex body of volume $1$ in $\mathbb R^n$ and let
$1\ls q, r\ls n$. There exists an absolute constant
$\overline{\beta }_2>0$ such that if $A$ is a subset of $K$
with $|A|\gr 1-e^{-\overline{\beta }_2q}$, then
\begin{equation}\label{eq:Lem2.3_1}
Z_p(K)\subseteq 2Z_p(\overline{A})
\end{equation}
for all $1\ls p\ls q$. Also, for the opposite inclusion, it
suffices to have $|A|\gr 2^{-\frac{r}{2}}$ to conclude that
\begin{equation}\label{eq:Lem2.3_2}
Z_p(\overline{A})\subseteq 2Z_p(K)
\end{equation}
for all $r\ls p\ls n$.
\end{lemma}

\noindent {\it Proof.} Let $\theta \in S^{n-1}$. Note that
\begin{equation}
h_{Z_p(\overline{A})}(\theta)=\left(\int_{\overline{A}}|\langle x,\theta\rangle |^pdx\right)^{1/p}=
\frac{1}{|A|^{\frac{1}{p}+\frac{1}{n}}}\left(\int_A|\langle x,\theta\rangle |^pdx\right)^{1/p}.
\end{equation}
We first prove $(\ref{eq:Lem2.3_2})$: since $A\subseteq K$ and
assuming that $|A|\gr 2^{-\frac{r}{2}}$, we have
\begin{align}
h_{Z_p(K)}(\theta) 
&= \left(\int_K|\langle x,\theta\rangle |^pdx\right)^{1/p} \gr\left(\int_A|\langle x,\theta\rangle |^pdx\right)^{1/p}
\\ \nonumber
&\gr 2^{-\frac{r}{2p}-\frac{r}{2n}}\left(\int_{\overline{A}}|\langle x,\theta\rangle|^pdx\right)^{1/p}\gr\frac{1}{2}h_{Z_p(\overline{A})}(\theta)
\end{align}
for all $r\ls p\ls n$. On the other hand, assuming that
$|A|\gr 1-e^{-\overline{\beta }_2q}$ and using the fact that
$\|\langle \cdot ,\theta\rangle\|_{2p}\ls 2\beta_2 \|\langle\cdot ,\theta\rangle \|_p$,
we have
\begin{align}
\int_K|\langle x, \theta \rangle |^pdx 
&= \int_A|\langle x, \theta \rangle |^pdx +\int_{K\setminus A}|\langle x, \theta\rangle |^pdx
\\ \nonumber
&\ls |A|^{1+\frac{p}{n}}\int_{\overline{A}}|\langle x, \theta \rangle |^pdx + |K\setminus A|^{1/2}\left(\int_K|\langle x,\theta\rangle |^{2p}dx\right)^{1/2}
\\ \nonumber
&\ls  \int_{\overline{A}}|\langle x, \theta \rangle |^pdx + e^{-\overline{\beta }_2q/2}(2\beta_2)^p\int_K|\langle x, \theta \rangle |^pdx
\\ \nonumber
&\ls \int_{\overline{A}}|\langle x, \theta \rangle |^pdx + \frac{1}{2}\int_K|\langle x,\theta\rangle |^pdx
\end{align}
for every $p\ls q$, if $\overline{\beta }_2>0$ is chosen large enough. This proves $(\ref{eq:Lem2.3_1})$. $\hfill\Box $

\section{Simple estimates for $I_1(K,Z_q^{\circ}(K))$}

In this section we give some upper and lower bounds for
$I_1(K,Z_q^{\circ}(K))$ which hold true for every isotropic convex
body $K$ in ${\mathbb R}^n$ and any $1\ls q\ls n$. In fact, our arguments
are quite direct and make use of estimates for simple parameters of
the bodies $Z_q(K)$, such as their radius or their volume, so that it
is straightforward to reach analogous upper and lower bounds for
$I_1(K,Z_q^{\circ}(M))$ in the more general case when $K$ and $M$
are not necessarily the same isotropic convex body.

\smallskip

Since $h_{Z_q(K)}(x)\ls R(Z_q(K))\| x\|_2$, we have that
\begin{equation}\label{eq:Sect3_1}
I_1(K,Z_q^{\circ}(K))\ls R(Z_q(K))\int_K\| x\|_2\,dx\ls R(Z_q(K))\sqrt{n}L_K,
\end{equation}
which, in combination with the fact that $R(Z_q(K))\ls \beta_1 qL_K$ (a direct
consequence of $(\ref{eq:Zq_1})$), leads to the bound
\begin{equation}\label{eq:Sect3_2}
I_1(K,Z_q^{\circ}(K))\ls \beta_1 q\sqrt{n}L_K^2.
\end{equation}
More generally, we have that
\begin{equation}\label{eq:Sect3_3}
I_1(K,Z_q^{\circ}(M))\ls R(Z_q(M))\int_K\| x\|_2\,dx\ls \beta_1 q\sqrt{n}L_K L_M.
\end{equation}
However, in the case that $M$ is an orthogonal transformation of $K$, the next lemma
shows that the average of the quantity $I_1(K,Z_q^{\circ}(M))$
can be bounded much more effectively than in $(\ref{eq:Sect3_3})$.

\begin{lemma}
Let $K$ be a centered convex body of volume $1$ in ${\mathbb R}^n$.
For every $2\ls q\ls n$,
\begin{equation}\left(\int_{O(n)}I_1^q\bigl(K,Z_q^{\circ}(U(K))\bigr)\,d\nu (U)\right)^{1/q}\ls C\sqrt{q/n}I_q^2(K),\end{equation}
where $C>0$ is an absolute constant.
\end{lemma}

\noindent {\it Proof.} We write
\begin{align}
\int_{O(n)}I_1^q\bigl(K,Z_q^{\circ }(U(K))\bigr)\,d\nu (U)
             &\ls \int_{O(n)}I_q^q\bigl(K,Z_q^{\circ }(U(K))\bigr)\,d\nu (U)
\\ \nonumber &= \int_{O(n)}\int_K\int_{U(K)}|\langle x,y\rangle |^qdy\,dx\,d\nu(U)
\\ \nonumber &= \int_K\int_K\int_{O(n)}|\langle x,Uy\rangle |^qd\nu (U)\,dy\,dx
\\ \nonumber &= \int_K\int_K\| y\|_2^q\int_{S^{n-1}}|\langle x,\theta \rangle|^qd\sigma (\theta )\,dy\,dx
\\ \nonumber &= c_{n,q}^q\int_K\int_K\| y\|_2^q\| x\|_2^qdy\,dx
\\ \nonumber &= c_{n,q}^qI_q^{2q}(K),
\end{align}
where $c_{n,q}\simeq\sqrt{q/n}$. $\hfill\Box $

\medskip

Recall that in the case that $K$ is isotropic, one has from \cite{Pa3} that
$I_q(K)\simeq\max\{\sqrt{n}L_K,R(Z_q(K))\}$. Then, Lemma 3.1
shows that, for every $2\ls q\ls n$,
\begin{equation}\left (\int_{O(n)}I_1^q\bigl(K,Z_q^{\circ }(U(K))\bigr)\,d\nu (U)\right )^{1/q}\ls C_1\max\{\sqrt{qn}, q^2\sqrt{q/n}\}\,L_K^2,\end{equation}
where $C_1>0$ is an absolute constant. Therefore, for every $2\ls q\ls\sqrt{n}$,
there exists a set $A_q\subseteq O(n)$ with $\nu (A_q)\gr 1-e^{-q}$ such that
$I_1\bigl(K,Z_q^{\circ }(U(K))\bigr)\ls C_2\sqrt{qn}\,L_K^2$ for all $U\in A_q$.
It is thus conceivable that there are properties of the bodies $Z_q(K)$ which
we can exploit to also bound $I_1(K,Z_q^{\circ }(K))$ more effectively than
in $(\ref{eq:Sect3_1})$ and $(\ref{eq:Sect3_2})$.

\medskip

We now pass to lower bounds; we will present three simple arguments. For the
first one we do not have to assume that $K$ or $M$ are in the isotropic position,
only that they are centered and have volume 1: from \cite[Corollary 2.2.a]{MP} we have that
\begin{equation}I_1(K,Z_q^{\circ}(M))= \int_K h_{Z_q(M)}(x)\,dx\gr \frac{n}{n+1}\, \frac{1}{|Z_q^{\circ}(M)|^{1/n}}.\end{equation}
Then, by the Blaschke--Santal\'{o} inequality, we get that
\begin{equation}I_1(K,Z_q^{\circ}(M))\gr c_1 n |Z_q(M)|^{1/n}\gr c_2\sqrt{qn}L_M\end{equation}
for all $2\ls q\ls \sqrt{n}$, because $|Z_q(M)|^{1/n}\gr c_3\sqrt{q/n}\,L_M$
for this range of values of $q$ by a recent result of Klartag and E. Milman (see \cite{KlEM}).
When $\sqrt{n}\ls q\ls n$, we have the weaker lower bound $|Z_q(M)|^{1/n}\gr c_4\sqrt{q/n}$,
which is due to Lutwak, Yang and Zhang (see \cite{LYZ1}). It follows
that $I_1(K,Z_q^{\circ }(M))\gr c_5\sqrt{qn}$ for this range of $q$.

For the second argument, we require that $K$ is isotropic and we write
\begin{align}
I_1(K,Z_q^{\circ}(M))
   &= \int_K h_{Z_q(M)}(x)\,dx=\int_K \max_{z\in Z_q(M)}|\langle x,z\rangle|\,dx
\\ \nonumber
   &\gr \max_{z\in Z_q(M)}\int_K|\langle x,z\rangle|\,dx \gr c\max_{z\in Z_q(M)}\|z\|_2L_K
\\ \nonumber
   &= c\, R(Z_q(M))L_K.
\end{align}

Finally, if $M$ is isotropic as well, we can use H\"{o}lder's inequality to get
\begin{align}
I_1(K,Z_q^{\circ}(M)) &= \int_K h_{Z_q(M)}(x)\,dx
\\ \nonumber
&\gr \int_K h_{Z_2(M)}(x)\,dx= \int_K \|x\|_2L_M\,dx\gr c\sqrt{n}L_K L_M.
\end{align}
All the estimates presented above are gathered in the next proposition.

\begin{proposition}
Let $K$ and $M$ be isotropic convex bodies in $\mathbb R^n$. For every $2\ls q\ls n$,
\begin{equation}
c_1\max\left\{\sqrt{n}L_K L_M, \sqrt{qn}, R(Z_q(M))L_K\right\}\ls
I_1(K,Z_q^{\circ}(M))\ls c_2q\sqrt{n}L_K L_M.
\end{equation}
In addition, if $2\ls q\ls \sqrt{n}$ then
\begin{equation}
c_1\max\left\{\sqrt{n}L_K L_M, \sqrt{qn}L_M\right\} \ls I_1(K,Z_q^{\circ}(M))\ls c_2q\sqrt{n}L_K L_M.
\end{equation}
\end{proposition}

The situation is more or less clear in the unconditional case.
Recall that a convex body $K$ in ${\mathbb R}^n$ is called
unconditional if it is symmetric with respect to all coordinate
hyperplanes (for some orthonormal basis of ${\mathbb R}^n$). Then,
it is easy to check that one can bring $K$ to the isotropic
position by applying an operator which is diagonal with respect to this
basis. It is also not hard to prove that the isotropic constant of
$K$ satisfies $L_K\simeq 1$. The upper bound follows from the
Loomis--Whitney inequality; see also \cite{BN1}. It is known (from
\cite{BN2}) that, for every $q\gr 2$, one has
$h_{Z_q(K)}(y)\ls c\sqrt{qn}\|y\|_{\infty }$,
where $c>0$ is an absolute constant. This leads us to the estimates
\begin{equation}
c_1\sqrt{qn}\ls I_1(K,Z_q^{\circ}(K))\ls c\sqrt{qn}\int_K\|x\|_{\infty}dx\ls c_2\sqrt{qn} \log\! n
\end{equation}
for all $2\ls q\ls n$ (the same estimates hold true for the quantity
$I_1(K,Z_q^{\circ}(M))$ when $M$ is too an unconditional isotropic convex body).

\section{The reduction}

Let $\kappa ,\tau >0$. Throughout this paper, we say that an isotropic convex body $K$ in
${\mathbb R}^n$ is $(\kappa ,\tau)$-regular if
\begin{equation}\label{eq:EntrEst}
\log N(K,tB_2^n)\ls \frac{\kappa n^2\log^2\! n}{t^2}\ \, \hbox{for all}\ t\gr \tau\sqrt{n\log\! n}.
\end{equation}
The purpose of this section is to present a reduction of the slicing
problem to the study of the quantity $I_1(K,Z_q^{\circ}(K))$ for $(\kappa ,\tau)$-regular
isotropic convex bodies: we show that any upper bound for $I_1(K,Z_q^{\circ}(K))$
immediately leads to an upper bound for the isotropic constant $L_K$ of a regular convex body $K$.
Note that the dependence seems to be nontrivial, in the sense that using the simple
estimates of Section 3 we can already recover the currently known bound
for $L_K$ with a loss of a logarithmic factor, while a small
(although not necessarily easy) improvement to those estimates
will also result in new bounds for $L_K$. In a sense, we will have
fully presented our reduction by the end of the next section, where
we provide a self-contained proof of the fact that there exist
regular isotropic convex bodies $K$ in $\mathbb R^n$ with $L_K\simeq L_n$.
First, let us see how the quantity $I_1(K,Z_q^{\circ}(K))$ and the isotropic
constant of a regular convex body $K$ are connected.

\begin{theorem}
There exists an absolute constant $\rho\in (0,1)$
with the following property: given $\kappa ,\tau\gr 1$, for every
$n\gr n_0(\tau )$ and every $(\kappa ,\tau)$-regular isotropic
convex body $K$ in ${\mathbb R}^n$ we have that, if $q\gr 2$ satisfies
\begin{equation}\label{eq:Th4.1_1}
2\ls q\ls \rho^2n\  {\textrm and}\ I_1(K,Z_q^{\circ }(K))\ls \rho nL_K^2,
\end{equation}
then
\begin{equation}
L_K^2\ls C\kappa\sqrt{\frac{n}{q}}\log^2\! n\,\max\left\{ 1,\frac{I_1(K,Z_q^{\circ}(K))}{\sqrt{qn}L_K^2}\right\}.
\end{equation}
\end{theorem}

\noindent {\it Proof.} We define a convex body $W$ in ${\mathbb R}^n$, setting
\begin{equation}W:=\{ x\in K: h_{Z_q(K)}(x)\ls C_1I_1(K,Z_q^{\circ }(K))\},\end{equation}
where $C_1=e^{2\overline{\beta}_2}$ and $\overline{\beta}_2>0$ is
the constant which appears in Lemma 2.3. From Markov's inequality
we have that $|W|\gr 1- e^{-2\overline{\beta}_2}$ and also
trivially that $|W|\gr 2^{-1}\gr 2^{-\frac{q}{2}}$ (as long as
$\overline{\beta}_2\gr 1$). Then we set
\begin{equation}K_1:=\overline{W}.\end{equation}
\noindent Applying both cases of Lemma 2.3 to the set $W$ with $p=2$, we see that
\begin{equation}\frac{1}{2}Z_2(K_1)\subseteq Z_2(K)\subseteq 2Z_2(K_1).\end{equation}
This implies that
\begin{equation}
\frac{1}{4}L_K^2=\frac{1}{4}\int_K\langle x,\theta\rangle^2dx\ls \int_{K_1}\langle x,\theta\rangle^2dx \ls 4\int_K\langle x,\theta\rangle^2dx=4L_K^2
\end{equation}
for every $\theta\in S^{n-1}$, and hence
\begin{equation}
\frac{nL_K^2}{4}\ls \sum_{i=1}^n\int_{K_1}\langle x,e_i\rangle^2dx=\int_{K_1}\|x\|_2^2dx\ls 4nL_K^2.
\end{equation}
We also have
\begin{equation}K_1=|W|^{-1/n}W\subseteq 2W\subseteq 2K,\end{equation}
thus for every $x\in K_1$ we have $x/2\in W$, and using $(\ref{eq:Lem2.3_2})$ of
Lemma 2.3 again, we can write
\begin{equation}\label{eqp:Th4.1_1}
h_{Z_q(K_1)}(x)\ls 2h_{Z_q(K)}(x)=4h_{Z_q(K)}(x/2)\ls 4C_1I_1(K,Z_q^{\circ }(K)).
\end{equation}
Finally,
\begin{equation}\label{eqp:Th4.1_2}
\log{N(K_1, tB_2^n)}\ls \log{N(2K,tB_2^n)}\ls \frac{4\kappa n^2\log^2\! n}{t^2},
\end{equation}
for all $t\gr 2\tau\sqrt{n\log\! n}$. We now write
\begin{equation}
nL_K^2 \ls 4\int_{K_1}\|x\|_2^2dx \ls 4\int_{K_1}\max_{z\in K_1}|\langle x,z \rangle |\,dx.
\end{equation}
$(\ref{eqp:Th4.1_2})$ tells us that for every $t\gr 2\tau\sqrt{n\log\! n}$, we can find $z_1,\ldots, z_{N_t}\in K_1$ such
that $K_1\subseteq \bigcup\limits_{i=1}^{N_t} (z_i+tB_2^n)$, and
$|N_t|\ls\exp\left (\frac{4\kappa n^2\log^2\! n}{t^2}\right )$. It follows that
\begin{equation}
\max_{z\in K_1}|\langle x,z\rangle |\ls\max_{1\ls i\ls N_t}|\langle x,z_i\rangle | + \max_{w\in tB_2^n} |\langle x,w\rangle|=\max_{1\ls i\ls N_t}|\langle x,z_i\rangle | + t\|x\|_2,
\end{equation}
and hence
\begin{align}\label{eqp:Th4.1_3}
nL_K^2 &\ls 4\int_{K_1}\max_{1\ls i\ls N_t}|\langle x,z_i\rangle |dx + 4t\int_{K_1}\|x\|_2dx
\\ \nonumber
       &\ls 4\int_{K_1}\max_{1\ls i\ls N_t}|\langle x,z_i\rangle |dx + 8t\sqrt{n}L_K.
\end{align}
We choose
\begin{equation}
t_0^2=16C_2\kappa\max\left\{ 1,\frac{I_1(K,Z_q^{\circ}(K))}{\sqrt{qn}L_K^2}\right\}\frac{n^{3/2}}{\sqrt{q}}\log^2\! n,
\end{equation}
where $C_2=16C_1\beta_2\overline{\beta}_1$ with $\beta_2$
the constant appearing in $(\ref{eq:Zq_2})$ and $\overline{\beta}_1$ the
constant from Lemma 2.1. With this choice of $t_0$, we have
\begin{equation}\label{eqp:Th4.1_4}
t_0^2\gr 16C_2\kappa\sqrt{\frac{n}{q}}\,n\log^2\! n\gr \frac{16C_2\kappa}{\rho}\,n\log^2\! n,
\end{equation}
as long as $q$ satisfies $(\ref{eq:Th4.1_1})$, and
\begin{equation}\label{eqp:Th4.1_5}
t_0^2\gr 16C_2\kappa\,\frac{I_1(K,Z_q^{\circ}(K))}{qL_K^2}\,n\log^2\! n.
\end{equation}
From $(\ref{eqp:Th4.1_4})$ it is clear that
\begin{equation}t_0^2\gr 16C_2\kappa \frac{n\log^2\! n}{\rho }\gr 4\tau^2n\log\! n,\end{equation}
provided that $n\gr n_0(\tau ,\kappa ,\rho)$, so the above argument,
leading up to $(\ref{eqp:Th4.1_3})$, remains valid for $t=t_0$. We also set
$p_0:=\frac{4\kappa n^2\log^2\! n}{t_0^2}$. Observe that $p_0\gr q$
(as long as $q$ is assumed to satisfy $(\ref{eq:Th4.1_1})$),
if $\rho$ is chosen properly: indeed, we have
$\max\left\{ 1,\frac{I_1(K,Z_q^{\circ}(K))}{\sqrt{qn}L_K^2}\right\}\ls \rho\sqrt{n/q}$,
and hence
\begin{equation}t_0^2\ls 16C_2\kappa \rho \frac{n^2\log^2\! n}{q}.\end{equation}
If we choose $\rho < 1/(4C_2)$, then we have
\begin{equation}
p_0 = \frac{4\kappa n^2\log^2\! n}{t_0^2}\gr \frac{4\kappa n^2q\log^2\! n}{16C_2\kappa\rho n^2\log^2\! n}
    = \frac{q}{4C_2\rho}\gr q.
\end{equation}
We therefore see that, using Lemma 2.1 with $q'=1$, we can write
\begin{equation}
\int_{K_1}\max_{1\ls i\ls N_{t_0}}|\langle x,z_i\rangle |dx\ls \overline{\beta}_1\max_{1\ls i\ls N_{t_0}}h_{Z_{p_0}(K_1)}(z_i)\ls \overline{\beta}_1\beta_2\frac{p_0}{q}\max_{1\ls i\ls N_{t_0}} h_{Z_q(K_1)}(z_i).
\end{equation}
Combining the above with $(\ref{eqp:Th4.1_3})$,
$(\ref{eqp:Th4.1_1})$ and the definition of $C_2$, we get
\begin{equation}\label{eqp:Th4.1_6}
nL_K^2\ls C_2\frac{p_0}{q}I_1(K,Z_q^{\circ}(K))+ 8t_0\sqrt{n}L_K. \end{equation}
Also, from $(\ref{eqp:Th4.1_5})$ and the definition of $p_0$, we have
\begin{equation}
C_2\frac{p_0}{q}I_1(K,Z_q^{\circ}(K))= \frac{4C_2\kappa I_1(K,Z_q^{\circ }(K))}{qt_0^2}\,n^2\log^2\! n\ls \frac{1}{4}\,nL_K^2.
\end{equation}
Therefore, $(\ref{eqp:Th4.1_6})$ becomes
\begin{equation}nL_K^2\ls C_3t_0\sqrt{n}L_K.\end{equation}
This gives us that
\begin{equation}
L_K^2\ls C_4\frac{t_0^2}{n}=C\kappa\max\left\{1,\frac{I_1(K,Z_q^{\circ}(K))}{\sqrt{qn}L_K^2}\right\}
\sqrt{\frac{n}{q}}\log^2\! n,
\end{equation}
as we desired. $\hfill\Box $

\section{Regular convex bodies with maximal isotropic constant}

Recall that $L_n:=\max\{L_K:K\;\hbox{is isotropic in}\,{\mathbb R}^n\}$.
In order to be able to use the argument of the previous
section to bound $L_n$, we need to establish the existence of $(\kappa ,\tau )$-regular convex bodies,
namely bodies satisfying $(\ref{eq:EntrEst})$, whose isotropic constant
is as ``close'' to $L_n$ as possible. The following theorem, formulated in the
more general setting of log-concave measures, was proven in $\cite{DP}$.

\begin{theorem}There exist absolute constants $\kappa ,\tau $ and $\delta
>0$ such that, for every $n\in {\mathbb N}$, there
exists an isotropic convex body $K$ in ${\mathbb R}^n$ with the
following properties:
\begin{enumerate}
\item[{\rm (i)}]
$L_K\gr \delta L_n$.
\item[{\rm (ii)}]
$\log{N(K, tB_2^n)}\ls\frac{\kappa n^2\log^2\! n}{t^2}$, for all $t\gr \tau\sqrt{n\log\! n}$.
\end{enumerate}
\end{theorem}

For the reader's convenience, we will give an outline of the proof
in the setting of convex bodies. First, we recall the following
theorem by Pisier which will be used in several steps of the
argument (see \cite{Pi} for a proof in the symmetric case; this
can easily be extended to the general case):

\begin{theorem}
Let $K$ be a centered convex body of volume $1$ in ${\mathbb R}^n$.
For every $\alpha\in (0,2)$ there exists an ellipsoid ${\cal{E}_{\alpha }}$
with $|{\cal{E}_{\alpha }}|=1 $ such that, for every $t\gr 1$,
\begin{equation}\label{eq:ThPisier}
\log N(K,t {\cal{E}_{\alpha}})\ls\frac{\kappa(\alpha)}{t^{\alpha}}n,
\end{equation}
where $\kappa (\alpha )>0$ is a constant depending only on $\alpha$.
\begin{remark}
{\rm One can take $\kappa(\alpha) \ls \frac{\kappa_1}{2-\alpha}$,
where $\kappa_1 >0$ is an absolute constant. An ellipsoid ${\cal E}_{\alpha }$
which satisfies $(\ref{eq:ThPisier})$ is called an $\alpha $-regular $M$-ellipsoid for $K$.}
\end{remark}
\end{theorem}

Secondly, let us gather some useful facts about ellipsoids in
${\mathbb R}^n$ that we are going to need for the proof of Theorem 5.1
(proofs for these facts can be found in \cite{DP}, \cite{KM1} and \cite{Z}).

\begin{lemma}
Let $\cal{E}$ be an ellipsoid in ${\mathbb R}^n$, then ${\cal E}
= T(B_2^n)$ for some $T\in GL(n)$. We denote the eigenvalues of
the matrix $\sqrt{T^{\ast}T}$ by $\lambda_1 \gr \cdots \gr \lambda_n>0$
$($recall that $T^{\ast}T$ is a symmetric, positive definite matrix$)$.
Then, for all $1\ls k\ls n-1$,
\begin{equation}\label{eq:LemEllips_1}
\max_{F \in G_{n,k}} |{\cal E} \cap F| = \max_{F \in G_{n,k}} |P_F({\cal E})|= \omega_k \prod_{i=1}^{k} \lambda_i
\end{equation}
and
\begin{equation}\label{eq:LemEllips_2}
\min_{F \in G_{n,k}} |{\cal E} \cap F| =\min_{F \in G_{n,k}} |P_F({\cal E})| =\omega_k  \prod_{i=n-k+1}^{n} \lambda_i.
\end{equation}
Also, if the dimension $n$ is even, we can find a subspace $F\in G_{n,n/2}$
such that $P_F({\cal E}) = \lambda_{n/2}B_F \ ( = \lambda_{n/2}B_2^n\cap F)$.
\end{lemma}

\smallskip

In view of the last part of Lemma 5.4, we choose to restrict ourselves to the
cases that the dimension $n$ is even, $n=2m$ for some $m\gr 1$, and prove Theorem
5.1 for those. However, as we will see in Remark 5.7, it is not hard to then
extend the theorem to all dimensions.

\smallskip

\noindent {\bf Proof of Theorem 5.1}. We start with an isotropic convex
body $K_1$ with $L_{K_1}\gr\delta_1 L_{2m}$, where $\delta_1\in (0,1)$. Then, one
has the following upper bound for the volume of sections of $K_1$.

\begin{lemma}
For every $k$-codimensional subspace $E$ of ${\mathbb R}^{2m}$,
$|K_1\cap E|^{1/k}\ls c_1(\delta_1 )$, where $c_1(\delta_1)>0$ depends only
on $\delta_1$.
\end{lemma}

\noindent {\it Proof.} Let $E$ be a $k$-codimensional subspace of
$\mathbb R^{2m}$, and denote its orthogonal subspace by $F$. We
consider the body $B_{k+1}(K_1,F)$, a convex body in the subspace
$F$ defined as in Subsection 2.2, and we recall that
\begin{equation}\label{eqp:Lem5.5_1}
c_1\frac{L_{B_{k+1}(K_1,F)}}{L_{K_1}}\ls |K_1\cap E|^{1/k}\ls c_2\frac{L_{B_{k+1}(K_1,F)}}{L_{K_1}}
\end{equation}
for some absolute constants $c_1, c_2$ independent of $m$ or $k$.
On the other hand, it is not hard to check that if
$k\ls j$ then $L_k\ls c_3L_j$ (see e.g. \cite[Theorem 4.2.2]{Gian}). Thus,
\begin{equation}L_{B_{k+1}(K_1,F)}\ls L_k\ls c_3L_{2m}=(c_3/\delta_1)L_{K_1},\end{equation}
and the lemma follows with $c_1(\delta_1 )=c_2c_3/\delta_1$. $\hfill\Box $

\medskip

We will now invoke Pisier's theorem to also give a lower bound
for the volume of $m$-dimensional sections of $K_1$ that contain
its barycenter.

\begin{lemma}For every $F\in G_{2m,m}$ we have $|K_1\cap F|^{1/m} \gr c_2(\delta_1)$,
where $c_2(\delta_1)>0$ depends only on $\delta_1$.
\end{lemma}

\noindent {\it Proof}. We consider an $\alpha $-regular
$M$--ellipsoid ${\cal E}_{\alpha }$ for $K_1$ (for the proof of this
lemma we could have fixed $\alpha =1$; however, some steps of this
more general argument will be needed again later). Set $t_{\alpha }=
\max\{1, [\kappa (\alpha)]^{1/\alpha }\}$. Then,
\begin{equation}\label{eqp:Lem5.6_1}
|P_F(K_1)|\ls N(K_1,t_{\alpha }{\cal E}_{\alpha })|P_F(t_{\alpha }{\cal E}_{\alpha })| \ls e^{2m}|P_F(t_{\alpha }{\cal E}_{\alpha})|
\end{equation}
for every $F\in G_{2m,m}$. We also need the Rogers-Shephard inequality (see \cite{RS2})
for both $K_1$ and ${\cal E}_{\alpha }$: since $|K_1|=|{\cal E}_{\alpha }|=1$, we know that
\begin{equation}\label{eqp:Lem5.6_2}
1=c_1\ls |K_1\cap F|^{1/m}|P_{F^{\perp }}(K_1)|^{1/m} \ls c_2,
\end{equation}
and similar estimates hold true for ${\cal E}_{\alpha }$
(see \cite{Spingarn} or \cite{MilPaj2} for the left hand side inequality).
The idea of the argument is the following:
inequality $(\ref{eqp:Lem5.6_2})$ helps us relate the volume of $m$-dimensional
sections of $K_1$ (or ${\cal E}_{\alpha}$) to that of $m$-dimensional
projections of $K_1$ (or ${\cal E}_{\alpha}$ respectively); an upper bound
for the former will give us a lower bound for the latter and vice versa.
Also, inequality $(\ref{eqp:Lem5.6_1})$ allows us to compare the maximum (or minimum) volume of the
$m$-dimensional projections of $K_1$ to the maximum (or minimum) volume of the
corresponding projections of ${\cal E}_{\alpha}$. However, as we recalled in Lemma 5.4,
the maximum volume of the $m$-dimensional projections of an ellipsoid
is the same as the maximum volume of its $m$-dimensional sections, so
we can use inequalities $(\ref{eqp:Lem5.6_1})$ and $(\ref{eqp:Lem5.6_2})$
once more to get from upper bounds for the volume of sections of $K_1$ to lower bounds.

\smallskip

We now give the precise argument: combining $(\ref{eqp:Lem5.6_2})$ with the conclusion
of Lemma 5.5, we see that $\min\limits_{F\in G_{2m,m}}|P_{F^{\perp}}(K_1)|^{1/m}\gr c_3(\delta_1)$.
We then get from $(\ref{eqp:Lem5.6_1})$ that $\min\limits_{F\in G_{2m,m}}|P_{F^{\perp
}}(t_{\alpha }{\cal E}_{\alpha})|^{1/m} \gr c_4(\delta_1)$.
Now, using $(\ref{eqp:Lem5.6_2})$ for ${\cal E}_{\alpha}$ we get $|{\cal E}_{\alpha}\cap
F|^{1/m}\ls c_5(\delta_1)t_{\alpha }$ for every $F\in G_{2m,m}$. But from $(\ref{eq:LemEllips_1})$ we have that
\begin{equation}\label{eqp:Lem5.6_3}
\max_{F\in G_{2m,m}}|P_F({\cal E}_{\alpha})|^{1/m} = \max_{F\in G_{2m,m}}|{\cal E}_{\alpha}\cap F|^{1/m}\ls c_5(\delta_1)t_{\alpha}.
\end{equation}
Using $(\ref{eqp:Lem5.6_1})$ once again, we get
$|P_F(K_1)|^{1/m}\ls c_6(\delta_1)t_{\alpha}^2$ for every $F\in
G_{2m,m}$. Inserting this estimate into $(\ref{eqp:Lem5.6_2})$, we see that
$|K_1\cap F|^{1/m}\gr c_7(\delta_1)/t_{\alpha}^2$ for every $F\in G_{2m,m}$. We
may choose $\alpha=1$ now, and complete the proof with
$c_2(\delta_1)=c_7(\delta_1)/t_1^2$. $\hfill\Box $

\medskip

\noindent {\bf Conclusion of the proof of Theorem 5.1}.
Let $\alpha\in (1,2)$ and let ${\cal E}_{\alpha}$ be an
$\alpha $-regular $M$--ellipsoid for $K$. Recall that
$|{\cal E}_{\alpha}|= 1$. Also, if  ${\cal E}_{\alpha} = T(B_2^n)= T(B_2^{2m})$,
let $\lambda_1 \gr \cdots \gr \lambda_{2m}>0$ be the eigenvalues of the
matrix $\sqrt{T^{\ast}T}$; observe from Lemma 5.4 that
\begin{equation}
|B_2^{m}|\prod_{i=m+1}^{2m}\lambda_i=\min_{F\in G_{2m,m}}|P_F({\cal E}_{\alpha })| \ls \max_{F\in G_{2m,m}}|P_F({\cal E}_{\alpha })| = |B_2^{m}|\prod_{i=1}^{m} \lambda_i.
\end{equation}
Using $(\ref{eqp:Lem5.6_1})$ and the conclusion of Lemma 5.6, we get
\begin{align}
|B_2^m|^{1/m}\lambda_m
  &\gr \min_{F\in G_{2m,m}}|P_F({\cal E}_{\alpha })|^{1/m}\gr
\frac{e^{-2}}{t_{\alpha}}\min_{F\in G_{2m,m}}|P_F(K_1)|^{1/m}
\\ \nonumber
  &\gr \frac{e^{-2}}{t_{\alpha }}\min_{F\in G_{2m,m}}|K_1\cap F|^{1/m}\gr \frac{c_8(\delta_1)}{t_{\alpha }},
\end{align}
and hence
\begin{equation}\lambda_m \gr \frac{c_9(\delta_1)}{t_{\alpha}}\sqrt{n}.\end{equation}
In a similar way, using $(\ref{eqp:Lem5.6_3})$, we see that
$|B_2^m|^{1/m}\lambda_m\ls \max\limits_{F\in G_{2m,m}}|P_F({\cal E}_{\alpha })|^{1/m}\ls c_5(\delta_1)t_{\alpha}$,
and hence $\lambda_m\ls c_{10}(\delta_1)t_{\alpha }\sqrt{n}$.
But from the last part of Lemma 5.4 we know
that there exists a subspace $F_0\in G_{2m,m}$ such that
$P_{F_0}({\cal{E}_{\alpha}})= \lambda_m B_{F_0}$, therefore,
\begin{equation}\label{eqp:Th5.1_1}
\frac{c_9(\delta_1)}{t_{\alpha }}\sqrt{n}B_{F_0} \subseteq P_{F_{0}}({\cal E}_{\alpha}) \subseteq c_{10}(\delta_1)t_{\alpha}\sqrt{n}B_{F_0}.
\end{equation}

\smallskip

Let $W:=\overline{B}_{m+1}(K_1,F_{0})$ and $K:= W\times U(W)$, where
$U\in O(2m)$ satisfies $U(F_0)=F_0^{\perp}$. Since $W$ is almost isotropic
and $L_{U(W)}=L_W$, from \cite[Lemma 1.6.6]{Gian} we see that
$K=W\times U(W)$ is an almost isotropic convex body in ${\mathbb R}^n \equiv {\mathbb R}^{2m}$
with $L_K=L_W$. We will show that $K$ satisfies (i) and (ii); the same conclusion
will then immediately follow (perhaps with slightly different constants for property (ii))
for any isotropic linear image $T(K)$ of $K$ satisfying $T(K)\simeq K$.

\medskip

\noindent {\it Proof of} (i): Since $L_K=L_W$, from $(\ref{eqp:Lem5.5_1})$ we get
\begin{equation}
L_K=L_W\gr c_2^{-1}L_{K_1} |K_1\cap F_{0}^{\perp}|^{1/m} \gr c_2^{-1}c_2(\delta_1)L_{K_1} \gr \delta L_n,
\end{equation}
where $\delta =\delta_1 c_2(\delta_1)/c_2$. For the last two inequalities we
have used Lemma 5.6 and the fact that $L_{K_1}\gr \delta_1 L_n$.

\smallskip

\noindent {\it Proof of} (ii): Using the fact that
$N(A\times A, B\times B)\ls N(A,B)^2$ for any two
nonempty sets $A,B$, and also the fact that
$B_2^m\times B_2^m\subseteq \sqrt{2} B_2^{2m}$,
we may write for any $s>0$,
\begin{equation}
N\bigl(K, s\sqrt{2n}B_2^n\bigr)\ls N\bigl(W\times U(W),s\sqrt{n}(B_{F_0}\times B_{F_0^{\perp}})\bigr)\ls N\bigl(W, s\sqrt{n}B_{F_0}\bigr)^2.
\end{equation}
From $(\ref{eq:KB2})$ we know that
\begin{equation}Z_m(\overline{B}_{m+1}(K_1,F_0)) \simeq |K_1\cap F_0^{\perp}|^{1/m} P_{F_0} (Z_m(K_1)),\end{equation}
therefore, using Lemmas 5.5, 5.6 and the fact that ${\rm conv}(C,-C)\simeq Z_m(C)$
for every centered convex body $C$ of volume $1$ in $F_0$ or in ${\mathbb R}^n$, we get
\begin{align}
{\rm conv}(W,-W)
    &\simeq Z_{m}(\overline{B}_{m+1}(K_1,F_0)) \simeq |K_1\cap F_0^{\perp}|^{1/m} P_{F_0} (Z_m(K_1))
\\ \nonumber
    &\simeq_{\delta_1} P_{F_0}({\rm conv}(K_1,-K_1)).
\end{align}
But then, recalling also $(\ref{eqp:Th5.1_1})$, we have for every $r>0$,
\begin{align}
N(W, c_{10}(\delta_1)t_{\alpha}r\sqrt{n}B_{F_0})
     &\ls  N\bigl({\rm conv}(W,-W),c_{10}(\delta_1)t_{\alpha}r\sqrt{n}B_{F_0}\bigr)
\\ \nonumber
     &\ls N\bigl({\rm conv}(W,-W),rP_{F_0}({\cal E}_{\alpha})\bigr)
\\ \nonumber
     &\ls N\bigl(c_{11}(\delta_1)P_{F_0}({\rm conv}(K_1,-K_1)),rP_{F_0}({\cal E}_{\alpha })\bigr)
\\ \nonumber
     &\ls N\bigl(c_{11}(\delta_1){\rm conv}(K_1,-K_1),r{\cal E}_{\alpha}\bigr)
\\ \nonumber
     &\ls N\bigl(K_1-K_1, c_{12}(\delta_1)r({\cal E}_{\alpha}-{\cal E}_{\alpha})\bigr)
\\ \nonumber
     &\ls N(K_1,c_{13}(\delta_1)r{\cal E}_{\alpha})^2
\end{align}
(note that for the last two inequalities we have also used that ${\cal E}_{\alpha}$ is
convex and symmetric, so ${\cal E}_{\alpha}-{\cal E}_{\alpha}= 2{\cal E}_{\alpha}$,
that $K_1$ is convex and contains the origin, so ${\rm conv}(K_1,-K_1)\subset K_1-K_1$,
as well as the fact that $N(A-A, B-B)\ls N(A,B)^2$).
It follows that
\begin{equation}\label{eqp:Th5.1_2}
N(K,t\sqrt{n}\,B_2^n)\ls N\left(K_1,\frac{c_{13}(\delta_1)t}{\sqrt{2}c_{10}(\delta_1)t_{\alpha}}\, {{\cal E}_{\alpha }}\right)^4
\end{equation}
for every $t>0$. Since ${{\cal E}_{\alpha}}$ is an $\alpha$-regular $M$--ellipsoid for $K_1$,
it remains to consider large enough $t\gr \tau (\delta_1,\alpha)$, where
\begin{equation}
\tau (\delta_1,\alpha) := \sqrt{2}c_{10}(\delta_1)t_{\alpha}/c_{13}(\delta_1) = t_{\alpha}/c_{14}(\delta_1),
\end{equation}
to deduce from $(\ref{eq:ThPisier})$ and $(\ref{eqp:Th5.1_2})$ that
\begin{equation}
\log N(K,t\sqrt{n}\,B_2^n)\ls 4\log N\left(K_1,\frac{c_{14}(\delta_1)t}{t_{\alpha}}\,{\cal E}_{\alpha}\right)\ls \frac{4\kappa(\alpha)t_{\alpha}^{\alpha}}{c_{14}^{\alpha}(\delta_1)}\,\!\frac{n}{t^{\alpha}}.
\end{equation}
Choosing $\alpha =2-\frac{1}{\log\! n}$, we have $\kappa(\alpha)\ls\kappa_1\log\! n$
and $t^{\alpha}\simeq t^2$ as long as, say, $t\ls n^2$. This completes the proof.
$\hfill\Box $

\smallskip

\begin{remark}
{\rm Now that we have proven the existence of an isotropic body $K$ in ${\mathbb R}^{2m}$
which has properties (i) and (ii) of Theorem 5.1, we can easily prove the existence of
such bodies in ${\mathbb R}^{2m-1}$ as well: just note that for every subspace $F\in G_{2m,2m-1}$
we have that $2L_K\ls |K\cap F^{\perp}|\ls 2R(K)$. Combining this with the properties $(\ref{eq:KB1}), (\ref{eq:KB2})$
for the almost isotropic convex body $\overline{B}_{2m}(K,F)$ in the $(2m-1)$--$\,\!$dimensional subspace $F$, we get that
\begin{equation}\label{eq:Rem2.7_1}
L_{\overline{B}_{2m}(K,F)}\simeq |K\cap F^{\perp}|^{\frac{1}{2m-1}}L_K \simeq L_K \gr \delta L_{2m}\gr c\delta L_{2m-1},
\end{equation}
and also that
\begin{equation}\label{eq:Rem2.7_2}
\overline{B}_{2m}(K,F)\simeq Z_{2m-1}\bigl(\overline{B}_{2m}(K,F)\bigr)\simeq
|K\cap F^{\perp}|^{\frac{1}{2m-1}}P_F(Z_{2m-1}(K))\simeq P_F(K).
\end{equation}
Since for every $t>0$, $N(P_F(K), tB_F)=N\bigl(P_F(K), tP_F(B_2^{2m})\bigr)\ls N(K, tB_2^{2m})$, we conclude that the body
$\overline{B}_{2m}(K,F)$ will also satisfy properties (i) and (ii) of Theorem 5.1 with perhaps slightly different,
but still independent of the dimension, constants $\kappa, \tau$ and $\delta$.}
\end{remark}

\medskip

In the statement of Theorem 5.1, we can add one more property about the radius of the
body $K$ that we look for: we can require that $R(K)\ls \gamma\sqrt{n}L_K$ where $\gamma >0$
is an absolute constant. The first step towards this is to use Bourgain's argument \cite{Bou}
which reduces the slicing problem to the study of bodies with small diameter; one can
prove the following fact (see e.g. \cite[Proposition 2.3.1]{Gian}).

\begin{lemma}
There exists an isotropic convex body $K_1$ in ${\mathbb R}^n$
with $L_{K_1}\gr\delta_1 L_n$ and $R(K_1)\ls \gamma_1\sqrt{n}L_{K_1}$,
where $\delta_1,\gamma_1>0$ are absolute constants.
\end{lemma}

Then, we can repeat the proof of Theorem 5.1 starting with the body
$K_1\subset {\mathbb R}^n = {\mathbb R}^{2m}$ given by Lemma 5.8. One has now that
$R(W)\ls c(\delta_1)\gamma_1\sqrt{n}L_{K_1}$: to see this, write
\begin{align}
R(W)
   &= R \bigl(\overline{B}_{m+1}(K_1,F_{0})\bigr)\ls c_1|K_1\cap F_0^{\perp}|^{1/m} R\bigl(P_{F_0}\bigl(Z_m(K_1)\bigr)\bigr)
\\ \nonumber
   &\ls c_2(\delta_1)R\bigl({\rm conv}(K_1,-K_1)\bigr) = c_2(\delta_1) R(K_1) \ls c_2(\delta_1)\gamma_1\sqrt{n}L_{K_1}.
\end{align}
It is also easy to check that $R(K)= R(W\times U(W))\simeq R(W)$,
hence $R(K)\ls \gamma\sqrt{n}L_K$ for some absolute constant
$\gamma>0$. Similarly for the odd dimensions, we see that for every $F\in G_{2m,2m-1}$,
\begin{equation}
R(\overline{B}_{2m}(K,F))\simeq R(P_F(K))\ls R(K)\ls c\gamma\sqrt{2m-1}\, L_{\overline{B}_{2m}(K,F)},
\end{equation}
where we have made use of $(\ref{eq:Rem2.7_1}), (\ref{eq:Rem2.7_2})$.
Thus, we can state the following version of Theorem 5.1.

\begin{theorem}
There exist absolute constants $\kappa ,\tau, \gamma $ and $\delta >0$
such that, for every $n\in {\mathbb N}$, there
can be found an isotropic convex body $K$ in $\mathbb R^n$ with
$R(K)\ls\gamma\sqrt{n}L_K$, $L_K\gr \delta L_n$, and the property that
\begin{equation}
\log{N(K, tB_2^n)}\ls\frac{\kappa n^2\log^2\! n}{t^2}\  \hbox{for all}\  t\gr \tau\sqrt{n\log\! n}.
\end{equation}
\end{theorem}

\begin{definition}{\rm Let ${\cal IK}(\kappa ,\tau ,\gamma ,\delta )$
denote the class of isotropic convex bodies whose existence is
established in Theorem 5.9. Let $\rho >0$ be the absolute constant
in Theorem 4.1. Then, we define $A(n,\kappa ,\tau ,\gamma ,\delta )$
to be the set of all $q\in [2,\rho^2n]$ for which there exists $K\in
{\cal IK}(\kappa ,\tau ,\gamma ,\delta )$ such that
$I_1(K,Z_q^{\circ }(K))\ls\rho nL_K^2$. Observe that already, by $(\ref{eq:Sect3_2})$,
$A(n,\kappa ,\tau ,\gamma ,\delta )$ can be shown to contain an interval
of the form $[2,c\sqrt{n}]$ where $c>0$ is an absolute constant. Clearly,
any improvement to the upper bound in $(\ref{eq:Sect3_2})$ will automatically give us that
$A(n,\kappa ,\tau ,\gamma ,\delta )$ contains an even larger part of $[2,\rho^2n]$.
For those $q$ we set
\begin{equation}
B(q)=\inf\left\{ \frac{I_1(K,Z_q^{\circ}(K))}{\sqrt{qn}L_K^2}:K\in {\cal IK}(\kappa ,\tau ,\gamma ,\delta)\right\}.
\end{equation}
Then, Theorem 4.1 implies the following:
for every $q\in A(n,\kappa ,\tau ,\gamma ,\delta )$,
\begin{equation}\delta^2L_n^2\ls C\kappa\sqrt{n/q}\log^2\! n\max\{ 1,B(q)\}.\end{equation}
In other words, we have:}
\end{definition}

\begin{theorem}
There exist absolute constants $\kappa ,\tau, \gamma $ and $\delta >0$ such that, for every $n\in {\mathbb N}$,
\begin{equation}
L_n^2\ls \min\left\{ \frac{C\kappa}{\delta^2}\sqrt{n/q}\log^2\! n\max\{1,B(q)\}:q\in A(n,\kappa ,\tau ,\gamma ,\delta )\right\}.
\end{equation}
\end{theorem}

The estimate $L_n\ls c\sqrt[4]{n}\log\! n$ is a direct consequence of
Theorem 5.11: observe that $B(2)\simeq 1$.

\section{Isotropic convex bodies with small diameter}

In \cite[Section 3]{GPV1} it is proven that for every isotropic
convex body $K$ there exists a second isotropic convex body $C$ with
bounded isotropic constant and the ``same behaviour" as $K$ with respect to
linear functionals.

\begin{theorem}
Let $K$ be an isotropic convex body in ${\mathbb R}^n$.
There exists an isotropic convex body $C$ in
${\mathbb R}^n$ with the following properties:
\begin{enumerate}
\item[{\rm (i)}]
$L_C\ls c_1$.
\item[{\rm (ii)}]
$c_2Z_q(C)\subseteq \frac{Z_q(K)}{L_K}+\sqrt{q}B_2^n\subseteq c_3Z_q(C)$ for all $1\ls q \ls n$.
\item[{\rm (iii)}]
$c_4I_q(C,W)\ls \frac{I_q(K,W)}{L_K}+I_q(D_n,W)\ls c_5I_q(C,W)$ for all $1\ls q\ls n$ and every symmetric convex body $W$ in ${\mathbb R}^n$.
\end{enumerate}
The constants $c_i$, $i=1,\ldots ,5$ are absolute positive constants.
\end{theorem}

The body $C$ is defined as the ``convolution" of $K$ with a multiple
of $B_2^n$. If we also assume that $K$ is symmetric, then using the
fact that $Z_n(C)\simeq C$ and $Z_n(K)\simeq K$, we see that
\begin{equation}\label{eq:Conv}
C\simeq\frac{K}{L_K}+D_n.
\end{equation}
From the previous section, we know that for our purposes
it is enough to study the quantity $I_1(K,Z_q^{\circ }(K))$
in the cases that $K$ is an isotropic symmetric convex body
with small diameter; that is, we can assume that $R(K)\ls \gamma\sqrt{n}L_K$
for some $\gamma\simeq 1$. The next proposition, which makes use
of Theorem 6.1, shows us that it even suffices to consider
isotropic convex bodies which are $c(\gamma)$-isomorphic to a
ball.

\begin{proposition}
Let $K$ be an isotropic symmetric convex body in ${\mathbb R}^n$
with $R(K)\ls \gamma\sqrt{n}L_K$. Then, there exists an isotropic
symmetric convex body $C$ such that:
\begin{enumerate}
\item[{\rm (i)}]
$L_C\ls c_6$,
\item[{\rm (ii)}]
$c_7D_n \subseteq C \subseteq c_8\gamma D_n$, and
\item[{\rm (iii)}]
$I_1(K,Z_q^{\circ }(K)) \ls c_9I_1(C,Z_q^{\circ }(C))L_K^2$ for all $1\ls q\ls n$,
\end{enumerate}
where $c_6,c_7,c_8, c_9>0$ are absolute constants.
\end{proposition}

\noindent {\it Proof.} We will use the fact that $w_q(Z_q(K))\simeq \sqrt{q/n}I_q(K)$, and hence
\begin{equation}c_1\sqrt{q}L_K\ls w_q(Z_q(K))\ls\gamma \sqrt{q}L_K \end{equation}
for all $1\ls q\ls n$. We consider the body $C$ defined by Theorem 6.1. It is clear
that $L_C\ls c_6$ for some absolute constant $c_6>0$. Since
$\frac{1}{L_K}Z_q(K)\subseteq c_3Z_q(C)$, we have
$c_3L_KZ_q^{\circ}(K)\supseteq Z_q^{\circ}(C)$, and hence
\begin{equation}
I_1(C, Z_q^{\circ}(C))\gr I_1(C, c_3L_KZ_q^{\circ}(K))= \frac{1}{c_3L_K}I_1(C, Z_q^{\circ}(K)).
\end{equation}
Applying the inequality $c_5I_1(C,W)\gr\frac{I_1(K,W)}{L_K}$ with
$W=Z_q^{\circ}(K)$, we get
\begin{equation}I_1(C, Z_q^{\circ}(C))\gr c_9\frac{I_1(K, Z_q^{\circ}(K))}{L_K^2},\end{equation}
with $c_9=(c_3c_5)^{-1}$ Finally, from
$(\ref{eq:Conv})$ and the fact that $\frac{K}{L_K}\subseteq\gamma \sqrt{n}B_2^n$,
we see that $c_7D_n\subseteq C\subseteq c_8\gamma D_n$. $\hfill\Box $

\medskip

In view of this result, we can give one more version of the
``reduction theorem" of Section 4.

\begin{definition}{\rm Let ${\cal IK}_{{\rm sd}}(\gamma )$
denote the class of isotropic convex bodies that satisfy
\begin{enumerate}
\item[{\rm (i)}]
$L_C\ls c_6$ and
\item[{\rm (ii)}]
$c_7D_n \subseteq C \subseteq c_8\gamma D_n$,
\end{enumerate}
where $c_i>0$ are absolute constants (e.g. the ones in Proposition
6.2). For every $2\ls q\ls n$, set
\begin{equation}
\Gamma (q)=\sup\left\{ \frac{I_1(K,Z_{q}^{\circ}(K))}{\sqrt{qn}}:K\in {\cal IK}_{{\rm sd}}(\gamma )\right\}.
\end{equation}
Then, Theorem 5.11 and Proposition 6.2 imply the following:}
\end{definition}

\begin{theorem}
There exist absolute constants $\kappa ,\tau, \gamma $ and $\delta >0$
such that, for every $n\in {\mathbb N}$,
\begin{equation}
L_n^2\ls \min\left\{ \frac{C\kappa}{\delta^2}\sqrt{n/q}\log^2\! n\, \Gamma (q):q\in A(n,\kappa ,\tau ,\gamma ,\delta)\right\}.
\end{equation}
\end{theorem}

In other words, studying the behaviour of
$\frac{I_1(K,Z_{q}^{\circ}(K))}{\sqrt{qn}}$ within the class
${\cal IK}_{{\rm sd}}(\gamma)$ is enough in order to understand
the behaviour of the parameter $B(q)$ of Section 5 as well as whether that
behaviour can lead to improved upper bounds for $L_n$.

\bigskip

\medskip

\noindent {\bf Acknowledgements}. The second named author is
partially supported by an NSF grant (DMS-0906150). The third named
author is supported by a scholarship of the University of Athens.
Part of this work was carried out while the second and the third
named authors were visiting the Isaac Newton Institute for
Mathematical Sciences in Cambridge; they would like to thank the
Institute for the hospitality and the organisers of the
{\it Discrete Analysis} Programme for their invitation.


\bigskip

\bigskip

\noindent \textsc{Apostolos \ Giannopoulos}: Department of
Mathematics, University of Athens, Panepistimioupolis 157-84,
Athens, Greece.

\smallskip

\noindent \textit{E-mail:} \texttt{apgiannop@math.uoa.gr}

\bigskip

\noindent \textsc{Grigoris \ Paouris}: Department of Mathematics,
Texas A \& M University, College Station, TX 77843 U.S.A.

\smallskip

\noindent \textit{E-mail:} {\tt grigoris\_paouris@yahoo.co.uk}

\bigskip

\noindent \textsc{Beatrice-Helen \ Vritsiou}: Department of
Mathematics, University of Athens, Panepistimioupolis 157-84,
Athens, Greece.

\smallskip

\noindent \textit{E-mail:} \texttt{bevritsi@math.uoa.gr}


\footnotesize
\begin{thebibliography}{00}
\footnotesize
\bibitem{Ball} {\rm K.\ M.\ Ball}, {\sl Logarithmically concave
functions and sections of convex sets in ${\mathbb R}^n$}, Studia
Math. {\bf 88} (1988), 69--84.

\bibitem{BN1} {\rm S.\ G.\ Bobkov and F.\ L.\ Nazarov}, {\sl On convex bodies and
log-concave probability measures with unconditional basis}, Geom.
Aspects of Funct. Analysis (Milman-Schechtman eds.), Lecture Notes
in Math. {\bf 1807} (2003), 53--69.
\bibitem{BN2} {\rm S.\ G.\ Bobkov and F.\ L.\ Nazarov},  {\sl Large deviations of
typical linear functionals on a convex body with unconditional
basis}, Stochastic Inequalities and Applications, Progr.\ Probab.\
56, Birkhauser, Basel, 2003, 3--13.

\bibitem{Bou} {\rm J. Bourgain}, {\sl On the distribution of polynomials on high dimensional
convex sets}, Geom. Aspects of Funct. Analysis (Lindenstrauss-Milman
eds.), Lecture Notes in Math. {\bf 1469} (1991), 127--137.
\bibitem{DP} {\rm N.\ Dafnis, G.\ Paouris}, {\sl Small ball probability estimates,
$\psi_2$-behavior and the hyperplane conjecture}, Journal of
Functional Analysis {\bf 258} (2010), 1933--1964.
\bibitem{Gian} {\rm A.\,Giannopoulos},
{\sl Notes on isotropic convex bodies}, {Warsaw University Notes}\,
(2003).
\bibitem{GPV1} {\rm A.\ Giannopoulos, G.\ Paouris and P.\ Valettas}, {\sl
On the existence of subgaussian directions for log-concave
measures}, Contemporary Mathematics {\bf 545} (2011), 103--122.
\bibitem{Kl} {\rm B.\ Klartag}, {\sl On convex perturbations with a bounded isotropic
constant}, Geom. Funct. Anal. {\bf 16} (2006), 1274--1290.
\bibitem{KlEM} {\rm B. Klartag and E. Milman}, {\sl Centroid Bodies and the
Logarithmic Laplace Transform - A Unified Approach},
arXiv:1103.2985v1.
\bibitem{KM1} {\rm B. Klartag and V. D. Milman},
{\sl Rapid Steiner Symmetrization of most of a convex body and the
slicing problem}, Combin. Probab. Comput. {\bf 14}, no. 5-6 (2005) 829--843.
\bibitem{LMS} {\rm A. Litvak, V. D. Milman and
G. Schechtman}, {\sl Averages of norms and quasi-norms}, Math. Ann.
{\bf 312} (1998), 95--124.
\bibitem{LYZ1} {\rm E. Lutwak, D. Yang and G. Zhang}, {\sl $L^p$
affine isoperimetric inequalities}, J. Differential Geom. {\bf 56}
(2000), 111--132.
\bibitem{MP} {\rm V.\ D.\ Milman and A.\ Pajor}, {\sl
Isotropic position and inertia ellipsoids and zonoids of the unit
ball of a normed $n$-dimensional space}, Geom. Aspects of Funct.
Analysis (Lindenstrauss-Milman eds.), Lecture Notes in Math. {\bf
1376} (1989), 64--104.
\bibitem{MilPaj2} {\rm V.D. Milman, A. Pajor}, {\sl Entropy and
Asymptotic Geometry of Non-Symmetric Convex Bodies}, Advances in
Mathematics,  {\bf 152} (2000), 314--335.
\bibitem{MS} {\rm V.D. Milman and G. Schechtman},
{\sl Asymptotic Theory of Finite Dimensional Normed Spaces}, Lecture
Notes in Math. {\bf 1200} (1986), Springer, Berlin.
\bibitem{Pa1} {\rm  G.\ Paouris}, {\sl $\psi_2$-estimates for linear functionals on zonoids},
Lecture Notes in Mathematics {\bf 1807} (2003), 211--222.
\bibitem{Pa3} {\rm  G.\ Paouris}, {\sl Concentration of mass in convex bodies},
Geometric and Functional Analysis  {\bf16} (2006), 1021--1049.
\bibitem{Pa4} {\rm  G.\ Paouris}, {\sl Small ball probability estimates for log-concave measures},
Trans. Amer. Math. Soc. (2011),
DOI:10.1090/S0002-9947-2011-05411-5.
\bibitem{Pi} {\rm G. Pisier}, {\sl The Volume
of Convex Bodies and Banach Space Geometry}, Cambridge Tracts in
Mathematics {\bf 94} (1989).
\bibitem{RS2} {\rm C. A. Rogers and G. C. Shephard}, {\sl Convex bodies associated with a
given convex body}, J. London Soc. {\bf 33} (1958), 270--281.
\bibitem{Spingarn} {\rm J. Spingarn}, {\sl An inequality for sections and projections of a convex set},
{Proc. Amer. Math. Soc.} {\bf 118} (1993), 1219--1224.
\bibitem{Z} {\rm C. Zong}, {\sl Strange phenomena in convex and discrete geometry},  Universitext,
Springer (2003).
\end{thebibliography}
\end{document}